\input amstex
\magnification=1200
\documentstyle{amsppt}
\NoRunningHeads
\NoBlackBoxes
\topmatter
\title Interactive games and representation theory. II.
A second quantization
\endtitle
\author\bf D.V.Juriev
\endauthor
\date math.RT/9808098\enddate
\endtopmatter
\document

This short note is devoted to a second quantization of a classical picture 
described in [1].

\subhead 1.1. Interactive games and intention fields\endsubhead

\definition{Definition 1 [1]} {\it An interactive system\/} (with $n$
{\it interactive controls\/}) is a control system with $n$ independent controls 
coupled with unknown or incompletely known feedbacks (the feedbacks, which are
called the {\it behavioral reactions}, as well as their couplings with 
controls are of a so complicated nature that their can not be described 
completely). {\it An interactive game\/} is a game with interactive controls 
of each player.
\enddefinition

Below we shall consider only deterministic and differential interactive
systems. For symplicity we suppose that $n=2$. In this case the general
interactive system may be written in the form:
$$\dot\varphi=\Phi(\varphi,u_1,u_2),\tag1$$
where $\varphi$ characterizes the state of the system and $u_i$ are
the interactive controls:
$$u_i(t)=u_i(u_i^\circ(t),\left.[\varphi(\tau)]\right|_{\tau\le t}),$$
i.e. the independent controls $u_i^\circ(t)$ coupled with the feedbacks on
$\left.[\varphi(\tau)]\right|_{\tau\le t}$. One may suppose that the
feedbacks are integrodifferential on $t$.

\proclaim{Proposition [1]} Each interactive system (1) may be transformed 
to the form (2) below (which is not, however, unique):
$$\dot\varphi=\tilde\Phi(\varphi,\xi),\tag2$$
where the magnitude $\xi$ (with infinite degrees of freedom as a rule) 
obeys the equation
$$\dot\xi=\Xi(\xi,\varphi,\tilde u_1,\tilde u_2),\tag3$$
where $\tilde u_i$ are the interactive controls of the form $\tilde 
u_i(t)=\tilde u_i(u_i^\circ(t); \varphi(t),\xi(t))$ (i.e. the feedbacks 
are on $\xi(t)$ as well as on $\varphi(t)$ and are differential on $t$).
\endproclaim

\remark{Remark 1} One may exclude $\varphi(t)$ from the feedbacks in
the interactive controls $\tilde u_i(t)$.
\endremark

\definition{Definition 2 [1]} The magnitude $\xi$ with its dynamical equations
(3) and its cont\-ri\-bution into the interactive controls $\tilde u_i$ will 
be called {\it the intention field}.
\enddefinition

\remark{Remark 2} The theorem holds true for the interactive games.
\endremark

\remark{Remark 3} In practice, the intention fields may be often considered 
as a field--theoretic description of subconscious individual and collective 
behavioral reactions. The intention fields realize a ``virtual dynamical 
memory'' for the transmission of cognitive data in the scheme of the 
accelerated nonverbal cognitive computer and telecommunications [2].
\endremark

\subhead 1.2. Second quantization of intention fields and inverse problem
of representation theory\endsubhead

Starting the classical picture sketched above one is able to perform its
second quantization. It means that the intention fields are quantized.
The quantum picture is more realistic for a description of processes of
the information transmission by intention fields in interactive games 
(cf.[3]).

The second quantization is deeply related to some inverse problems of the
representation theory.

\definition{Definition 3 [1]} {\it The main inverse problem of representation 
theory for the interactive system (1)\/} (or for the interactive game) is
\roster
\item"(1)" to write the system (1) in the form (2);
\item"(2)" to determine the geometrical and algebraical structure of
the intention field;
\item"(3)" to find the algebraic structure, which ``governs'' the dynamics
(3).
\endroster
\enddefinition

\remark{Remark 4} The solution of the main inverse problem of representation 
theory for the interactive system may use {\sl a posteriori} data on the 
system.
\endremark

So to perform a second quantization of the intention field it is sufficient
to solve the dynamical inverse problem of representation theory [4] for
the interactive system (game) (1). 

\remark{Remark 5} This article as well as [1] may be regarded as a 
realization of some ideas of [5] on the mathematics beyond the conventional 
one based on the well-known forms of perception such as vision. Thus, the 
theory of interactive games plays a role for the subconscious individual
and collective behavioral reactions analogous to geometry for a visual 
perception.
\endremark

\remark{Remark 6} I suppose that it is reasonable to describe `Chi' (`Tsi')
of Chineese tradition and its manifestations as (perhaps, quantum) intention 
fields (cf.[3]), at least partially. Apparently, in concrete situations
we should have deal with complex couplings of intention fields with
physical magnitudes (i.e. the complexes of $\xi$ and $\phi$).
\endremark 

\Refs
\roster
\item" [1]" Juriev D., Interactive games and representation theory.
E-print: math.FA/9803020.
\item" [2]" Juriev D., Droems: experimental mathematics, informatics and
infinite dimensional geometry. Report RCMPI-96/05+ [e-version: 
cs.HC/9809119].
\item" [3]" Juriev D.V., Complex projective geometry and quantum projective
field theory. Theor. Math.Phys. 101(3) (1994) 1387-1403.
\item" [4]" Juriev D., Dynamical inverse problem of representation theory
and noncommutative geometry [in Russian]. Fundam.Prikl.Matem. 4(1) (1998)
[e-version: funct-an/9507001].
\item" [5]" Saaty T.L., Speculating on the future of mathematics. 
Appl.Math.Lett. 1 (1988) 79-82.
\endroster
\endRefs
\enddocument